\documentclass[12pt]{article}
\usepackage{graphicx, amssymb}

\newtheorem{lem}{Lemma}[section]
\newtheorem{thm}[lem]{Theorem}
\newtheorem{pro}[lem]{Proposition}
\newtheorem{cor}[lem]{Corollary}

\newcommand{\qed}{\hfill{$\Box$}}
\newcommand{\ZZ}{{\mathbb{Z}}}
\newcommand{\U}{{\textsf{U}}}
\renewcommand{\L}{{\textsf{L}}}
\newcommand{\D}{{\textsf{D}}}

\begin{document}

\baselineskip=12pt

\centerline{\Large\bf A Simple Proof of the Aztec Diamond Theorem}

\vspace{1.27cm} \centerline{{\sc Sen-Peng
Eu}$^{1,}$\footnote{Partially supported by National Science
Council, Taiwan (NSC 92-2119-M-390-001).} and {\sc Tung-Shan
Fu}$^{2,}$\footnote{Partially supported by National Science
Council, Taiwan (NSC 92-2115-M-251-001).}}

\vspace{0.6cm} \noindent{$^1${\footnotesize{\em Department of
Applied Mathematics, National University of Kaohsiung, Kaohsiung
811, Taiwan}}}

\noindent{\small{\tt  \mbox{ }speu@nuk.edu.tw}}

\vspace{0.3cm} \noindent{$^2${\footnotesize{\em Mathematics
Faculty, National Pingtung Institute of Commerce, Pingtung 900,
Taiwan}}}

\noindent{\small{\tt  \mbox{ }tsfu@npic.edu.tw}}

\vspace{0.85cm}
\begin{abstract} Based on a bijection between domino tilings of an Aztec
diamond and non-intersecting lattice paths, a simple proof of the
Aztec diamond theorem is given in terms of Hankel determinants of
the large and small Schr\"oder numbers.

\vspace{0.6cm}

\noindent{\em MSC2000: 05A15}

\vspace{0.2cm} \noindent{\em Keywords:} Aztec diamond, Hankel
matrix, Schr\"oder numbers, lattice paths
\end{abstract}

\vspace{0.6cm} \baselineskip=20pt
\section{Introduction}

\begin{figure}[h]
\begin{center}
\includegraphics[width=4in]{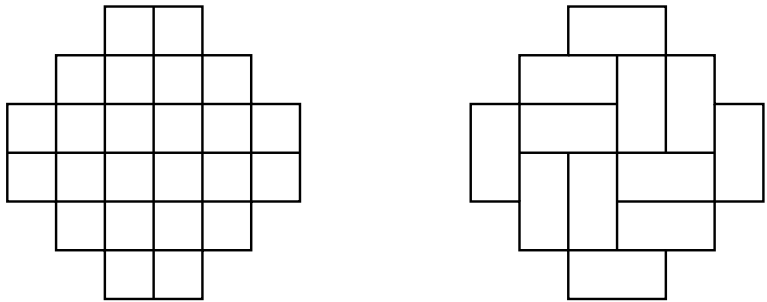}
\end{center}
\caption{the Az($3$) and a domino tiling} \label{fig:diamond}
\end{figure}

The {\em Aztec diamond} of order $n$, denoted by Az($n$), is
defined as the union of all the unit squares with integral corners
$(x,y)$ satisfying $|x|+|y|\le n+1$. A {\em domino} is simply a
1-by-2 or 2-by-1 rectangles with integral corners. A {\em domino
tiling} of a region $R$ is a set of non-overlapping dominos the
union of which is $R$. Figure \ref{fig:diamond} shows the Aztec
diamond of order 3 and a domino tiling. The Aztec diamond theorem,
which is first proved by Elkies {\em et al.} in \cite{EKLP1},
indicates that the number $a_n$ of domino tilings of the Aztec
diamond of order $n$ is $2^{n(n+1)/2}$. They gave four proofs by
relating the tilings to alternating sign matrices, monotone
triangles, representations of general linear groups, and domino
shuffling. Other approaches to this theorem appeared in \cite{BK,
Ciucu, Kuo}. Ciucu \cite{Ciucu} derived the recurrence relation
$a_n=2^na_{n-1}$ by means of perfect matchings of celluar graphs.
Kuo \cite{Kuo} developed a method, called graphical condensation,
to derive the recurrence relation $a_na_{n-2}=2a_{n-1}^2$, for
$n\ge 3$. Recently, Brualdi and Kirkland \cite{BK} gave a proof by
considering a matrix of order $n(n+1)$ the determinant of which
gives $a_n$. In this note we give a proof in terms of Hankel
determinants of the large and small Schr\"oder numbers based on a
bijection between the domino tilings of an Aztec diamond and
non-intersecting lattice paths.

Recall the {\em large} {\em Schr\"oder numbers} $\{r_n\}_{n\ge
0}:=\{1,2,6,22,90,394,1806,\dots\}$ and the {\em small Schr\"oder
numbers} $\{s_n\}_{n\ge 0}:=\{1,1,3,11,45,197,903,\dots\}$. Among
many other combinatorial structures, the $n$-th large Schr\"oder
number $r_n$ counts the number of lattice paths in the plane
$\ZZ\times\ZZ$ from $(0,0)$ to $(2n,0)$ using {\em up} steps
$(1,1)$, {\em down} steps $(1,-1)$, and {\em level} steps $(2,0)$
that never pass below the $x$-axis. Such a path is called a {\em
large Schr\"oder path} of length $n$ (or a {\em large} $n$-{\em
Schr\"oder path} for short). Let $\U$, $\D$, and $\L$ denote an
up, down, and level step, respectively.
 Note that the terms of
$\{r_n\}_{n\ge 1}$ are twice of those in $\{s_n\}_{n\ge 1}$.
Consequently, the $n$-th small Schr\"oder number $s_n$ counts the
number of large $n$-Schr\"oder paths without level steps on the
$x$-axis, for $n\ge 1$. Such a path is called a {\em small}
$n$-{\em Schr\"oder path}. Refer to \cite[Exercise 6.39]{Stanley}
for more information.

Our proof relies on the determinants of the following {\em Hankel
matrices} of the large and small Schr\"oder numbers
\[H^{(1)}_n:=\left[
           \begin{array}{cccc}
                r_1 &  r_2 & \cdots & r_n \\
                r_2 &  r_3 & \cdots & r_{n+1} \\
                \vdots & \vdots &   & \vdots \\
                r_n &  r_{n+1} & \cdots & r_{2n-1}
                \end{array}
           \right],
 \quad
 G^{(1)}_n:=\left[
           \begin{array}{cccc}
                s_1 &  s_2 & \cdots & s_n \\
                s_2 &  s_3 & \cdots & s_{n+1} \\
                \vdots & \vdots &   & \vdots \\
                s_n &  s_{n+1} & \cdots & s_{2n-1}
                \end{array}
           \right].
\]
Note that $H^{(1)}_n=2G^{(1)}_n$. Using a method of Gessel and
Viennot \cite{GV}, we associate the determinants of $H^{(1)}_n$
and $G^{(1)}_n$ with the numbers of $n$-tuples of non-intersecting
large and small Schr\"oder paths, respectively. How to derive the
determinants of $H^{(1)}_n$ and $G^{(1)}_n$ and how to establish
bijections between domino tilings of an Aztec diamond and
non-intersecting large Schr\"oder paths are given in the next
section.

\vspace{1cm}
\section{A proof of the Aztec diamond theorem}

Let $\Pi_n$ (resp. $\Omega_n$) denote the set of $n$-tuples
$(\pi_1,\dots,\pi_n)$ of large Schr\"oder paths (resp. small
Schr\"oder paths) satisfying the following two conditions.
\begin{enumerate}
\item[(A1)] The path $\pi_i$ goes from $(-2i+1,0)$ to $(2i-1,0)$, for $1\le i\le n$, and \\
\vspace{-1cm} \item[(A2)] any two paths $\pi_i$ and $\pi_j$ do not
intersect.
\end{enumerate}

\begin{figure}[h]
\begin{center}
\includegraphics[width=6.3in]{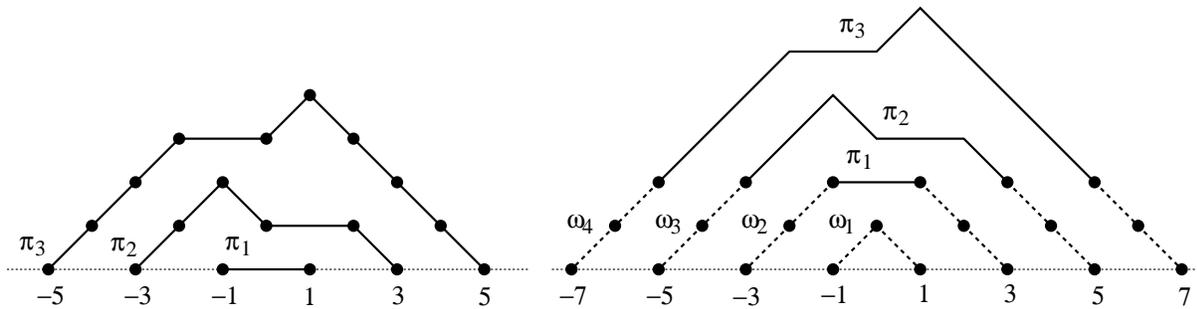}
\end{center}
\caption{\small a triple $(\pi_1,\pi_2,\pi_3)\in\Pi_3$ and the
corresponding quadruple
$(\omega_1,\omega_2,\omega_3,\omega_4)\in\Omega_4$}
\label{fig:omega}
\end{figure}

\bigskip
There is an immediate bijection $\phi$ between $\Pi_{n-1}$ and
$\Omega_n$, for $n\ge 2$, which carries
$(\pi_1,\dots,\pi_{n-1})\in\Pi_{n-1}$ into
$\phi((\pi_1,\dots,\pi_{n-1}))=(\omega_1,\dots,\omega_n)\in\Omega_n$,
where $\omega_1=\U\D$ and $\omega_i=\U\U\pi_{i-1}\D\D$ (i.e.,
$\omega_i$ is obtained from $\pi_{i-1}$ with 2 up steps attached
in the beginning and 2 down steps attached in the end, and then
rises above the $x$-axis), for $2\le i\le n$. For example, on the
left of Figure \ref{fig:omega} is a triple
$(\pi_1,\pi_2,\pi_3)\in\Pi_3$. The corresponding quadruple
$(\omega_1,\omega_2,\omega_3,\omega_4)\in\Omega_4$ is shown on the
right. Hence, for $n\ge 2$, we have
\begin{equation} \label{eqn:omega}
|\Pi_{n-1}|=|\Omega_n|.
\end{equation}

\bigskip
For a permutation $\sigma=z_1z_2\cdots z_n$ of $\{1,\dots,n\}$,
the {\em sign} of $\sigma$, denoted by sgn$(\sigma)$, is defined
by sgn$(\sigma):=(-1)^{{\rm inv}(\sigma)}$, where inv$(\sigma):=$
Card$\{(z_i,z_j)|\, i<j$ and $z_i>z_j\}$ is the number of {\em
inversions} of $\sigma$. Using the technique of a sign-reversing
involution over a signed set, we prove that the cardinalities of
$\Pi_n$ and $\Omega_n$ coincide with the determinants of
$H^{(1)}_n$ and $G^{(1)}_n$, respectively. Following the same
steps as \cite[Theorem 5.1]{DS}, a proof is given here for
completeness.

\bigskip
\begin{pro} \label{pro:Hankel-det} For $n\ge 1$, we have
\begin{enumerate}
\item $|\Pi_n|=\det(H^{(1)}_n)=2^{n(n+1)/2}$, and\item
$|\Omega_n|=\det(G^{(1)}_n)=2^{n(n-1)/2}$.
\end{enumerate}
\end{pro}

\medskip
\noindent{\em Proof:} For $1\le i\le n$, let $A_i$ denote the
point $(-2i+1,0)$ and let $B_i$ denote the point $(2i-1,0)$. Let
$h_{ij}$ denote the $(i,j)$-entry of $H^{(1)}_n$. Note that
$h_{ij}=r_{i+j-1}$ is equal to the number of large Schr\"oder
paths from $A_i$ to $B_j$. Let $P$ be the set of ordered pairs
$(\sigma, (\tau_1,\dots,\tau_n))$, where $\sigma$ is a permutation
of $\{1,\dots,n\}$, and  $(\tau_1,\dots,\tau_n)$ is an $n$-tuple
of large Schr\"oder paths such that $\tau_i$ goes from $A_i$ to
$B_{\sigma(i)}$. According to the sign of $\sigma$, the ordered
pairs in $P$ are partitioned into $P^+$ and $P^-$. Then
\[
\det(H^{(1)}_n)=\sum_{\sigma\in S_n} \mbox{\rm sgn}(\sigma)
\prod_{i=1}^n h_{i,\sigma(i)}=|P^+|-|P^-|.
\]
If there exists a sign-reversing involution $\varphi$ on $P$, then
$\det(H^{(1)}_n)$ is equal to the number of fixed points of
$\varphi$. Let $(\sigma,(\tau_1,\dots,\tau_n))\in P$ be such a
pair that at least two paths of $(\tau_1,\dots,\tau_n)$ intersect.
Choose the first pair $i<j$ in lexical order such that $\tau_i$
intersects $\tau_j$. Construct new paths $\tau_i'$ and $\tau_j'$
by switching the tails after the last point of intersection of
$\tau_i$ and $\tau_j$. Now $\tau_i'$ goes from $A_i$ to
$B_{\sigma(j)}$ and $\tau_j'$ goes from $A_j$ to $B_{\sigma(i)}$.
Since $\sigma\circ(ij)$ carries $i$ into $\sigma(j)$, $j$ into
$\sigma(i)$, and $k$ into $\sigma(k)$, for $k\neq i,j$, we define
\[\varphi((\sigma,(\tau_1,\dots,\tau_n)))=(\sigma\circ(ij),(\tau_1,\dots,\tau_i',\dots,\tau_j',\dots,\tau_n)).
\]
Clearly, $\varphi$ is sign-reversing. Since the first intersecting
pair $i<j$ is not affected by $\varphi$,  $\varphi$ is an
involution. The fixed points of $\varphi$ are the pairs
$(\sigma,(\tau_1,\dots,\tau_n))\in P$ such that $\sigma$ is the
identity, and  $\tau_1,\dots,\tau_n$ do not intersect, i.e.,
$(\tau_1,\dots,\tau_n)\in\Pi_n$. Hence $\det(H^{(1)}_n)=|\Pi_n|$.
By the same argument, we have $\det(G^{(1)}_n)=|\Omega_n|$. It
follows from (\ref{eqn:omega}) and the identity
$H^{(1)}_n=2G^{(1)}_n$ that
\[|\Pi_n|=\det(H^{(1)}_n)=2^n\cdot
\det(G^{(1)}_n)=2^n|\Omega_n|=2^n|\Pi_{n-1}|.\] Note that
$|\Pi_1|=2$, and hence, by induction, the assertions (i) and (ii)
follow. \qed

\bigskip
To prove the Aztec diamond theorem, we shall establish a bijection
between $\Pi_n$ and the set of domino tilings of Az($n$) based on
an idea, due to D. Randall, mentioned in \cite[Solution of
Exercise 6.49]{Stanley}.

\bigskip
\begin{pro} \label{pro:bijetion} There is a bijection between
the set of domino tilings of the Aztec diamond of order $n$ and
the set of $n$-tuples $(\pi_1,\dots,\pi_n)$ of large Schr\"oder
paths satisfying the conditions (A1) and (A2).
\end{pro}

\medskip
\noindent{\em Proof:} Given a tiling $T$ of Az($n$), we associate
$T$ with an $n$-tuple $(\tau_1,\dots,\tau_n)$ of non-intersecting
paths as follows. Let the rows of Az($n$) be indexed by
$1,2,\dots,2n$ from bottom to top.  For $1\le i\le n$, define a
path $\tau_i$ from the center of the left-hand edge of the $i$-th
row to the center of the right-hand edge of the $i$-th row.
Namely, each step of the path is from the center of a domino edge
(where a domino is regarded as having six edges of unit length) to
the center of another edge of the some domino $D$, such that the
step is symmetric with respect to the center of $D$. One can check
that for each tiling there is a unique such an $n$-tuple
$(\tau_1,\dots,\tau_n)$ of paths, moreover, any two paths
$\tau_i$, $\tau_j$ of which do not intersect. Conversely any such
$n$-tuple of paths corresponds to a unique domino tiling of
Az($n$)  (note that any domino not on these paths is horizontal).

  To establish a mapping $\psi$, for $1\le i\le n$, we
form a large Schr\"oder path $\pi_i$ from $\tau_i$ with $i-1$ up
steps attached in the beginning of $\tau_i$ and with $i-1$ down
steps attached in the end (and then raise $\pi_i$ above the
$x$-axis), and define $\psi(T)=(\pi_1,\dots,\pi_n)$. One can
verify that the $n$-tuple $(\pi_1,\dots,\pi_n)$ of large
Schr\"oder paths satisfies the conditions (A1) and (A2), and hence
$\psi(T)\in\Pi_n$. To find $\psi^{-1}$, we can retrieve an
$n$-tuple $(\tau_1,\dots,\tau_n)$ of non-intersecting paths, which
corresponds to a unique domino tiling of Az($n$), from each
$n$-tuple $(\pi_1,\dots,\pi_n)$ of large Schr\"oder paths
satisfying the conditions (A1) and (A2) by a reverse procedure.
\qed

\begin{figure}[ht]
\begin{center}
\includegraphics[width=5in]{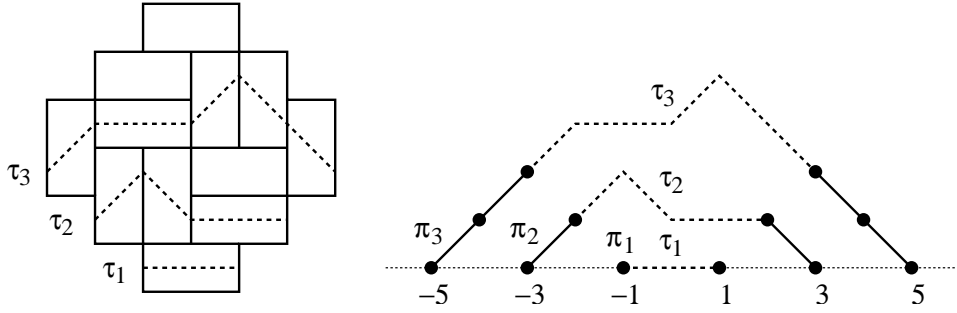}
\end{center}
\caption{\small a tiling of Az($3$) and the corresponding triple
of non-intersecting Schr\"oder paths} \label{fig:Aztec-Schroder}
\end{figure}

\bigskip
For example, on the left of Figure \ref{fig:Aztec-Schroder} is a
tiling $T$ of Az($3$) and the associated triple
$(\tau_1,\tau_2,\tau_3)$ of non-intersecting paths. On the right
of Figure \ref{fig:Aztec-Schroder} is the corresponding triple
$\psi(T)=(\pi_1,\pi_2,\pi_3)\in\Pi_3$ of large Schr\"oder paths.

\bigskip
By Propositions \ref{pro:Hankel-det} and \ref{pro:bijetion}, we
deduce the Aztec diamond theorem anew.

\bigskip
\begin{thm} \label{thm:Aztec-diamond}{\rm\bf (Aztec diamond theorem)} The number of
domino tilings of the Aztec diamond of order $n$ is
$2^{n(n+1)/2}$.
\end{thm}

\bigskip
\noindent{\bf Remark:} The proof of Proposition
\ref{pro:Hankel-det} relies on the recurrence relation
$\Pi_{n}=2^n\Pi_{n-1}$ essentially, which is derived by means of
the determinants of the Hankel matrices $H^{(1)}_n$ and
$G^{(1)}_n$. We are interested to hear a purely combinatorial
proof of this recurrence relation.

\bigskip
 In a similar manner we derive the
determinants of the Hankel matrices of large and small Schr\"oder
paths of the forms
\[H^{(0)}_n:=\left[
           \begin{array}{cccc}
                r_0 &  r_1 & \cdots & r_{n-1} \\
                r_1 &  r_2 & \cdots & r_n \\
                \vdots & \vdots &   & \vdots \\
                r_{n-1} &  r_n & \cdots & r_{2n-2}
                \end{array}
           \right],
 \quad
 G^{(0)}_n:=\left[
           \begin{array}{cccc}
                s_0 &  s_1 & \cdots & s_{n-1} \\
                s_1 &  s_2 & \cdots & s_n \\
                \vdots & \vdots &   & \vdots \\
                s_{n-1} &  s_n & \cdots & s_{2n-2}
                \end{array}
           \right].
\]

\bigskip
\begin{pro} For $n\ge 1$, $\det(H^{(0)}_n)=\det(G^{(0)}_n)=2^{n(n-1)/2}$.
\end{pro}

\medskip
\noindent{\em Proof:} Let $\Pi_n^*$ (resp. $\Omega_n^*$) be the
set of $n$-tuples $(\mu_0,\mu_1,\dots,\mu_{n-1})$ of large
Schr\"oder paths (resp. small Schr\"oder paths) satisfying the two
conditions (i) the path $\mu_i$ goes from $(-2i,0)$ to $(2i,0)$,
for $0\le i\le n-1$, and (ii) any two paths $\mu_i$ and $\mu_j$ do
not intersect. Note that $\mu_0$ degenerates into a single point
and that $\Pi_n^*$ and $\Omega_n^*$ are identical since for any
$(\mu_0,\mu_1,\dots,\mu_{n-1})\in\Pi_n^*$ all of the paths $\mu_i$
have no level steps on the $x$-axis. By a similar argument of
Proposition \ref{pro:Hankel-det}, we have
$\det(H^{(0)}_n)=|\Pi_n^*|=|\Omega_n^*|=\det(G^{(0)}_n)$.
Moreover, there is a bijection $\rho$ between $\Pi_{n-1}$ and
$\Pi_n^*$, for $n\ge 2$, which carries
$(\pi_1,\dots,\pi_{n-1})\in\Pi_{n-1}$ into
$\rho((\pi_1,\dots,\pi_{n-1}))=(\mu_0,\mu_1,\dots,\mu_{n-1})\in\Pi_n^*$,
where $\mu_0$ is the origin and $\mu_i=\U\pi_i\D$, for $1\le i\le
n-1$. The assertion follows from Proposition
\ref{pro:Hankel-det}(i). \qed

\begin{figure}[ht]
\begin{center}
\includegraphics[width=6in]{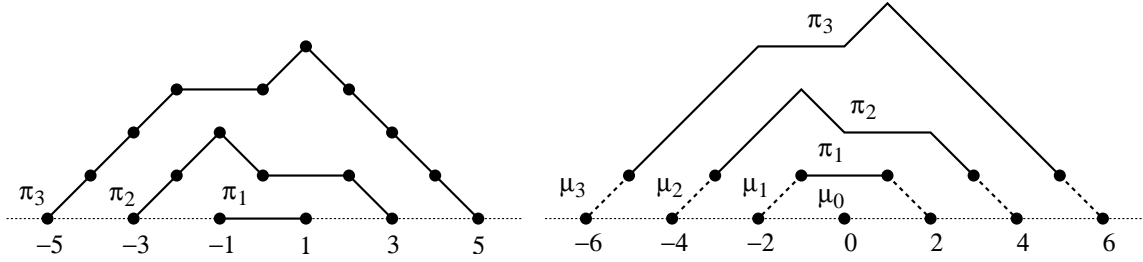}
\end{center}
\caption{\small a triple $(\pi_1,\pi_2,\pi_3)\in\Pi_3$ and the
corresponding quadruple $(\mu_0,\mu_1,\mu_2,\mu_3)\in\Pi_4^*$}
\label{fig:mu}
\end{figure}

\bigskip
For example, on the left of Figure \ref{fig:mu} is a triple
$(\pi_1,\pi_2,\pi_3)\in\Pi_3$ of non-intersecting large Schr\"oder
paths. The corresponding quadruple
$(\mu_0,\mu_1,\mu_2,\mu_3)\in\Pi_4^*$ is shown on the right.

\bigskip
Hankel matrices $H^{(0)}_n$ and $H^{(1)}_n$ may be associated with
any given sequence of real numbers. As noted by Aigner in
\cite{Aigner} that the sequence of determinants
\[\det(H^{(0)}_1),\det(H^{(1)}_1),\det(H^{(0)}_2),\det(H^{(1)}_2),\dots\]
uniquely determines the original number sequence provided that
$\det(H^{(0)}_n)\neq 0$ and $\det(H^{(1)}_n)\neq 0$, for all $n\ge
1$, we have a characterization of large and small Schr\"oder
numbers.

\bigskip
\begin{cor} The following results hold.
\begin{enumerate}
\item The large Schr\"oder numbers $\{r_n\}_{n\ge 0}$ are the
unique sequence with the Hankel determinants
$\det(H^{(0)}_n)=2^{n(n-1)/2}$ and $\det(H^{(1)}_n)=2^{n(n+1)/2}$,
for all $n\ge 1$.

\item The small Schr\"oder numbers $\{s_n\}_{n\ge 0}$ are the
unique sequence with the Hankel determinants
$\det(G^{(0)}_n)=\det(G^{(1)}_n)=2^{n(n-1)/2}$, for all $n\ge 1$.
\end{enumerate}
\end{cor}

\end{document}